\documentclass[11pt]{article}
\usepackage{enumerate}
\usepackage{amssymb,a4wide,latexsym,makeidx,epsfig,fleqn}
\usepackage{amsthm}
\usepackage{amsmath}
\usepackage{enumerate}
\usepackage{graphicx}
\usepackage{subfigure}
\usepackage{float}
\usepackage{caption}
\usepackage[colorlinks, linkcolor=black, anchorcolor=black, citecolor=blue]{hyperref}%
\newtheorem{theorem}{Theorem}[section]

\newtheorem{lemma}[theorem]{Lemma}

\newtheorem{corollary}[theorem]{Corollary}

\begin{document}
\textwidth 150mm \textheight 225mm
\title{Maxima of the $Q$-index of 2 leaves-free graphs with given size\thanks{Supported by the National Natural Science Foundation of China (No. 12271439).}}
\author{{Yuxiang Liu$^{a,b}$, Ligong Wang$^{a,b,}$\footnote{Corresponding author.}, Xiaolong Jia$^{c}$}\\
{\small $^{a}$  School of Mathematics and Statistics, Northwestern Polytechnical University,}\\{\small  Xi'an, Shanxi 710129, P.R. China.}\\ {\small $^{b}$ Xi'an-Budapest Joint Research Center for Combinatorics, Northwestern Polytechnical University,}\\{\small  Xi'an, Shanxi 710129, P.R. China.}\\ {\small $^{c}$ Xi'an Dilun Biotechnology Co., Ltd}\\{\small  Xi'an, Shanxi 710000, P.R. China.}\\
{\small E-mail: yxliumath@163.com, lgwangmath@163.com, xlxljj@126.com}}
\date{}
\maketitle
\begin{center}
\begin{minipage}{135mm}
\vskip 0.3cm
\begin{center}
{\small {\bf Abstract}}
\end{center}
{\small The $Q$-index of graph $G$ is the largest eigenvalue of the signless Laplacian matrix of $G$. Wang [Discrete Appl. Math. 356(2024)] proved the sharp upper bounds on the $Q$-index of leaf-free graphs with given size and characterized the corresponding extremal graphs. A graph is $2$ leaves-free if it has no two pendent vertices. In this paper, we give sharp upper bounds on the $Q$-index of 2 leaves-free graphs with given size and characterize the corresponding extremal graphs.
}
\end{minipage}
\end{center}
\section{Introduction}
Let $G$ be an undirected simple graph with vertex set $V(G)=\{v_{1},\cdots, v_{n}\}$ and edge set
$E(G)$, where $e(G)$ is the number of edges of $G$. For $v\in V(G)$, the neighborhood $N_{G}(v)$ of $v$
is $\{u: uv\in E(G)\}$ and the degree $d_{G}(v)$ of $v$ is $N_{G}(u)$. Denote by $\Delta(G)$ and $\delta(G)$ the maximum
and minimum degree of $G$, respectively. Denote by $P_{n}, C_{n}$, $S_{n}$ the path, the cycle and the star graph of order $n$, respectively. Denote by $D_{a,b}$ the double star graph consisting of connecting two center vertices of $K_{1,a}$ and $K_{1,b}$. For $A, B \subseteq V(G)$, $e(A)$ denotes the number of the edges of $G$ with both endvertices in $A$ and $e(A, B)$ denotes the number of the edges of $G$ with one endvertex in $A$ and the other in $B$. For two vertex disjoint graphs $G$ and $H$, we denote by $G$ and $H$ and $G\vee H$ the union of $G$ and $H$, and the join of $G$ and $H$, i.e., joining every vertex of $G$ to every vertex of $H$, respectively. Denote by $kG$ the union of $k$ disjoint copies of $G$. Let $S_{n}^{+}$ be the graph obtained from $K_{1,n-1}$ by adding one edge within its independent sets. For any nonisolated vertex $v$ of $G$, the average degree of $v$ is denoted by $m(v)=\frac{1}{d(v)}\sum_{u\in N(v)}d(u)$. For graph notation and terminology undefined here,we refer the readers to \cite{BoMu1}.

For a graph $G$, let $A(G)$ denote its adjacency matrix and $D(G)$ denote the diagonal matrix of vertex degrees. The matrix
$Q(G)=D(G) + A(G)$ is called the signless Laplacian matrix (or the $Q$-matrix) of $G$. The largest eigenvalues of $A(G)$ and $Q(G)$ are called the index (denoted by $\rho(G)$) and the $Q$-index (denoted by $q(G)$) of $G$, respectively. By the well-known Perron--Frobenius theorem, if $G$ is connected, then $Q(G)$ is irreducible and there exists a unique positive unit eigenvector $X = (x_{1}, x_{2}, \cdots, x_{n})^{T}$ corresponding to $q(G)$ such that $q(G)=X^{\prime}QX=\sum_{ij\in E(G)}(x_{i}+x_{j})^{2}$.

The investigation on the upper or lower bounds of the index and the $Q$-index of graphs is an important and classic
topic in the theory of graph spectra. For example, one may see \cite{FSi,LiLF,Ni5,Ni7,Ni3} and the references therein. Specially, the problem of characterizing the graph with maximal index with given size was posed by Brualdi and Hoffman \cite{BrH} as a conjecture, and solved by Rowlinson \cite{Row}. Motivated by this problem, a good deal of attention has been devoted to determining the graphs with maximum index ($Q$-index) among a given class of graphs with given size. See, for
example, \cite{GuoZ,LouG,ZhXL1,ZhXL2} and the references therein. Zhai, Xue and Luo \cite{ZhXL2} determined the maximal $Q$-index of graphs with given size and clique number (resp., chromatic number). Zhai, Xue and Liu \cite{ZhXL1} determined the maximal $Q$-index of graphs with given size and matching number. Lou, Guo and Wang \cite{LouG} characterized the unique graph of size $m$ with diameter at least $d$ having the largest $Q$-index and $L$-index. Jia, Li and Wang \cite{JiaLW} ordered all the graphs with given size $m$ and diameter $d$ from the second to the ($\lfloor\frac{d}{2}\rfloor+1$)th via their largest $Q$-indices. It is worth noting that all those extremal graphs shown in \cite{LouG,ZhXL1,ZhXL2} contain pendent vertices and recently Wang \cite{Wsj} determined the maximal $Q$-index of leaf-free graphs with given size and characterized the corresponding extremal graphs completely, which motivates us to consider those graphs with no 2 pendent vertices, that is, the 2 leaves-free graphs.

\noindent\begin{theorem}\label{th:ch-1.1.}{\rm(}$\cite{Wsj}${\rm)} Let $G$ be a leaf-free graph with size $m\geq12$.

{\rm(}$i${\rm)} If $m=3k$, then $q(G)\leq \frac{2k+3+\sqrt{4k^{2}-4k+9}}{2}$, with equality if and only if $G\cong K_{1}\vee kP_{2}$.

{\rm(}$ii${\rm)} If $m=3k+1$, then $q(G)\leq \alpha$, with equality if and only if $G\cong K_{1}\vee ((k-2)P_{2}\cup S_{4})$, where $\alpha$ is the largest root of the equation $x^{4}-(2k+9)x^{3}+(16k+23)x^{2}-(34k+19)x+20k-4=0$.

{\rm(}$iii${\rm)} If $m=3k+2$, then $q(G)\leq \beta$, with equality if and only if $G\cong K_{1}\vee ((k-1)P_{2}\cup S_{3})$,where $\beta$ is the largest root of the equation $x^{4}-(2k+9)x^{3}+(14k+26)x^{2}-28(k+1)x+16k+8=0$.
\end{theorem}

\noindent\begin{theorem}\label{th:ch-1.2.} Let $G$ be a 2 leaves-free graph with size $m\geq17$.

{\rm(}$i${\rm)} If $m=3k$, then $q(G)\leq \gamma$, with equality if and only if $G\cong K_{1}\vee ((k-2)P_{2}\cup S_{3}\cup P_{1})$, where $\gamma$ is the largest root of the equation $x^{5}-(2k+8)x^{4}+(14k+19)x^{3}-(28k+11)x^{2}+(16k-15)x+12=0$.

{\rm(}$ii${\rm)} If $m=3k+1$, then $q(G)\leq \frac{2k+3+\sqrt{4k^{2}+4k+9}}{2}$, with equality if and only if $G\cong K_{1}\vee (kP_{2}\cup P_{1})$.

{\rm(}$iii${\rm)} If $m=3k+2$, then $q(G)\leq \xi$, with equality if and only if $G\cong K_{1}\vee ((k-2)P_{2}\cup S_{4}\cup P_{1})$, where $\xi$ is the largest root of the equation $x^{5}-(2k+10)x^{4}+(16k+33)x^{3}-(36k+42)x^{2}+24kx+18=0$.
\end{theorem}

\section{Preliminary}
In this section, we introduce some lemmas and notations. For a graph $G$ and a subset $S\subseteq V(G)$, let $G[S]$ be the subgraph of $G$ induced by $S$. For two vertex subsets $S$ and $T$ of $V(G)$ (where $S\cap T$ may not be empty), let $e(S,T)$ denote the number of edges with one endpoint in $S$ and the other in $T$. $e(S,S)$ is simplified by $e(S)$. For a vertex $u\in V(G)$, let $N(u)$ be the neighborhood of the vertex $u$ of $G$, $N^{2}(u)$ be the set of vertices of distance two to $u$. Let $N_{S}(v)=N(v)\cap S$ and $d_{S}(v)=|N_{S}(v)|$. Let $N_{+}(u)=N(u)\setminus N_{0}(u)$ where $N_{0}(u)$ is the set of isolated vertices of $G[N(u)]$. Let $W=V(G)\setminus N[u]$. Let $X$ be the Perron vector of $G$ with coordinate $x_{v}$ corresponding to the vertex $v\in V(G)$ and $u^{\ast}$ be a vertex satisfying $x_{u^{\ast}}=\max\{x_{v}|v\in V (G)\}$.

\noindent\begin{lemma}\label{le:ch-2.1.}{\rm(}$\cite{Ste}${\rm)}
Let $H$ be a proper subgraph of a connected graph $G$. Then $q(H)<q(G)$. In particular, $q(G)\geq q(K_{1,\Delta(G)})=\Delta(G)+1$, with equality if and only if $G\cong K_{1,\Delta(G)}$.
\end{lemma}
\noindent\begin{lemma}\label{le:ch-2.2.} {\rm(}$\cite{HyZ}${\rm)} Let $G$ be a connected graph, $u$ and $v$ be two distinct vertices of $G$. Suppose that $v_{i}\in N(v)\setminus N(u)$ ($1\leq i\leq s$) and $X=(x_{1}, x_{2}, \ldots, x_{n})^{T}$ be the Perron vector of $Q(G)$. Let $G^{\prime}=G-\sum_{i=1}^{s}v_{i}v+\sum_{i=1}^{s}v_{i}u$. If $x_{u}\geq x_{v}$, then $q(G)<q(G^{\prime})$.
\end{lemma}

\noindent\begin{lemma}\label{le:ch-2.3.}{\rm(}$\cite{FeY}${\rm)} Let $G$ be a graph without isolated vertices. Then we have that $$q(G)\leq max_{v\in V(G)}\{d(v)+m(v)\},$$ with equality if and only if $G$ is regular or bi-regular.
\end{lemma}

\noindent\begin{lemma}\label{le:ch-2.4.}{\rm(}$\cite{Wsj}${\rm)} Let $G$ be a graph with size $m$. Then for any non-isolated  vertex $v\in V(G)$, $$m(v)\leq min\{\frac{2m}{d(v)}-1, \frac{m}{d(v)}+\frac{d(v)-1}{2}\}.$$
\end{lemma}

\section{Proof of Theorem 1.2.}
In this section, we give the  proof of Theorem 1.2.

\noindent\begin{lemma}\label{le:ch-2.9.}{\rm(}$\cite{Wsj}${\rm)}
Let $\mathcal{F}_{s,r}$ be the set of forest $F$ with $r$ components such that $\delta(F)=1$ and $|E(F)|=s$. Then for any $G\in\{K_{1}\vee F: F\in \mathcal{F}_{s,r}\}$, we have that $q(G)\leq q(K_{1}\vee (S_{s-r+2}\cup (r-1)P_{2}))$, with equality if and only if $G\cong K_{1}\vee (S_{s-r+2}\cup (r-1)P_{2})$.
\end{lemma}

\noindent\begin{corollary}\label{le:ch-2.10.}
Let $\mathcal{F}_{s,r}$ be the set of forest $F$ with $r$ components such that $\delta(F)=1$ and $|E(F)|=s$. Then for any $G\in\{K_{1}\vee (F\cup P_{1}): F\in \mathcal{F}_{s,r}\}$, we have that $q(G)\leq q(K_{1}\vee (S_{s-r+2}\cup (r-1)P_{2}\cup P_{1}))$, with equality if and only if $G\cong K_{1}\vee (S_{s-r+2}\cup (r-1)P_{2}\cup P_{1})$.
\end{corollary}

\noindent\begin{lemma}\label{le:ch-2.5.} Let $G$ be a 2 leaves-free graph of size $m$, then we have $\Delta(G)\leq \lfloor\frac{2m+1}{3}\rfloor$. Furthermore, the equality holds if and only if
\[
\begin{cases}
G\in \{K_{1}\vee (kP_{2}\cup P_{1}),K_{1}\vee ((k-2)P_{2}\cup S_{3}\cup P_{1})\}&\text{$m=3k$};\\
G\cong K_{1}\vee (kP_{2}\cup P_{1})&\text{$m=3k+1$};\\
G\in \{L_{1},L_{2}, K_{1}\vee ((k-1)P_{2}\cup S_{3}), K_{1}\vee ((k-2)P_{2}\cup S_{4}\cup P_{1}),K_{1}\vee ((k-2)P_{2}\cup P_{4}\cup P_{1}),&\\ K_{1}\vee ((k-3)P_{2}\cup 2S_{3}\cup P_{1})\}&\text{$m=3k+2$}.
\end{cases}
\]
\end{lemma}

\noindent{\bf Proof.} Since $G$ is 2 leaves-free graph, we obtain there exists at most a vertex $u$ satisfying $d(u)=1$ and $|V(G)|\geq \Delta(G)+1$. Furthermore, we have
$$2m=\sum_{v\in V(G)}d(v)\geq \Delta(G)+1+2|V(G)-2|=3\Delta(G)-1.$$ Hence, we have $\Delta(G)\leq \lfloor\frac{2m+1}{3}\rfloor$.

The sufficient condition can be proved by calculate, so we only prove the necessity.

(i) If $m=3k$ and $d(w)=\Delta(G)=\lfloor\frac{2m+1}{3}\rfloor=2k$, then we have that $|V(G-w)|\geq d(w)=2k$ and $|E(G-w)|=k$.

{\bf Case 1.} If $|V(G-w)|\geq 2k+1$, let $W=G\setminus N[w]$, then $|W|\geq1$ and we consider theses cases as follows.

{\bf Subcase 1.1.} If $\delta(G)\geq2$, then

$$2m=6k=\sum_{v\in V(G)}d(v)\geq 2k+2\times 2k+2|W|=6k+2|W|,$$ and hence $|W|=0$, a contradiction.

{\bf Subcase 1.2.} If there exists a vertex with degree 1, then
either $$2m=6k=\sum_{v\in V(G)}d(v)\geq 2k+2\times (2k-1)+1+2|W|=6k-1+2|W|,$$ or $$2m=6k=\sum_{v\in V(G)}d(v)\geq 2k+2\times 2k+1+2|W-1|=6k-1+2|W|, $$ and hence $|W|=0$, a contradiction.

{\bf Case 2.} If $|V(G-w)|=2k$, i.e., $G\cong K_{1}\vee (G-w)$, then we consider theses cases as follows.

{\bf Subcase 2.1.} If $\delta(G)\geq2$, then $$2m=6k=\sum_{v\in V(G)}d(v)\geq 2k+2\times 2k=6k, $$ and hence the degree sequence of $G-w$ is $\underbrace{(1,1,\cdots,1}_{2k})$. Hence $G-w\cong kP_{2}$ and $G\cong K_{1}\vee kP_{2}$.

{\bf Subcase 2.2.} If there exists a vertex with degree 1, then $$2m=6k=\sum_{v\in V(G)}d(v)\geq 2k+2\times (2k-1)+1=6k-1, $$ and hence there exists at most a vertex with degree 2 in $G-w$. Thus, we have the degree sequence of $G-w$ is $(2, \underbrace{1,1,\cdots,1}_{2k-2},0)$. Thus $G\cong K_{1}\vee ((k-2)P_{2}\cup S_{3}\cup P_{1})$.

(ii) If $m=3k+1$ and $d(w)=\Delta(G)=\lfloor\frac{2m+1}{3}\rfloor=2k+1$, then we have that $|V(G-w)|\geq d(w)=2k+1$ and $|E(G-w)|=k$.

{\bf Case 1.} If $|V(G-w)|\geq 2k+2$, then $|W|\geq1$ and we consider these cases as follows.

{\bf Subcase 1.1.} If $\delta(G)\geq2$, then $$2m=6k+2=\sum_{v\in V(G)}d(v)\geq 2k+1+2\times (2k+1)+2|W|=6k+3+2|W|,$$ and hence $|W|<0$, a contradiction.

{\bf Subcase 1.2.} If there exists a vertex with degree 1, then either $$2m=6k+2=\sum_{v\in V(G)}d(v)\geq 2k+1+2\times 2k+1+2|W|=6k+2+2|W|, $$ or $$2m=6k+2=\sum_{v\in V(G)}d(v)\geq 2k+1+2\times (2k+1)+1+2(|W|-1)=6k+2+2|W|, $$ and hence $|W|=0$, a contradiction.

{\bf Case 2.} If $|V(G-w)|=2k+1$, i.e., $G\cong K_{1}\vee (G-w)$, then we consider these cases as follows.

{\bf Subcase 2.1.} If $\delta(G)\geq2$, then $$2m=6k+2=\sum_{v\in V(G)}d(v)\geq 2k+1+2\times (2k+1)=6k+3, $$ a contradiction.

{\bf Subcase 2.2.} If there exists a vertex with degree 1, then $$2m=6k+2=\sum_{v\in V(G)}d(v)\geq 2k+1+2\times 2k+1=6k+2, $$ and hence the degree sequence of $G-w$ is $\underbrace{(1,1,\cdots,1}_{2k},0)$. Thus, $G\cong K_{1}\vee (kP_{2}\cup P_{1})$.

(iii) If $m=3k+2$ and $d(w)=\Delta(G)=\lfloor\frac{2m+1}{3}\rfloor=2k+1$, then we have that $|V(G-w)|\geq d(w)=2k+1$ and $|E(G-w)|=k+1$.

{\bf Case 1.} If $|V(G-w)|\geq 2k+2$, then $|W|\geq1$ and we consider theses cases as follows.

{\bf Subcase 1.1.} If $\delta(G)\geq2$, then $$2m=6k+4=\sum_{v\in V(G)}d(v)\geq 2k+1+2\times (2k+1)+2|W|=6k+3+2|W|,$$ and hence $|W|=0$, a contradiction.

{\bf Subcase 1.2.} If there exists a vertex with degree 1, then either $$2m=6k+4=\sum_{v\in V(G)}d(v)\geq 2k+1+2\times 2k+1+2|W|=6k+2+2|W|, $$ or $$2m=6k+4=\sum_{v\in V(G)}d(v)\geq 2k+1+2\times (2k+1)+1+2(|W|-1)=6k+2+2|W|, $$ and hence $|W|=1$.
Furthermore, the degree sequence of $G-w$ is $(\underbrace{1,1,\cdots,1}_{2k},0,2)$ or $(\underbrace{1,1,\cdots,1}_{2k+1},1)$. Thus, $G\in \{L_{1},L_{2}\}$ (see Fig. 1).
\begin{figure}[H]
\begin{centering}
 \subfigure[$L_{1}$]{
  \includegraphics[scale=0.25]{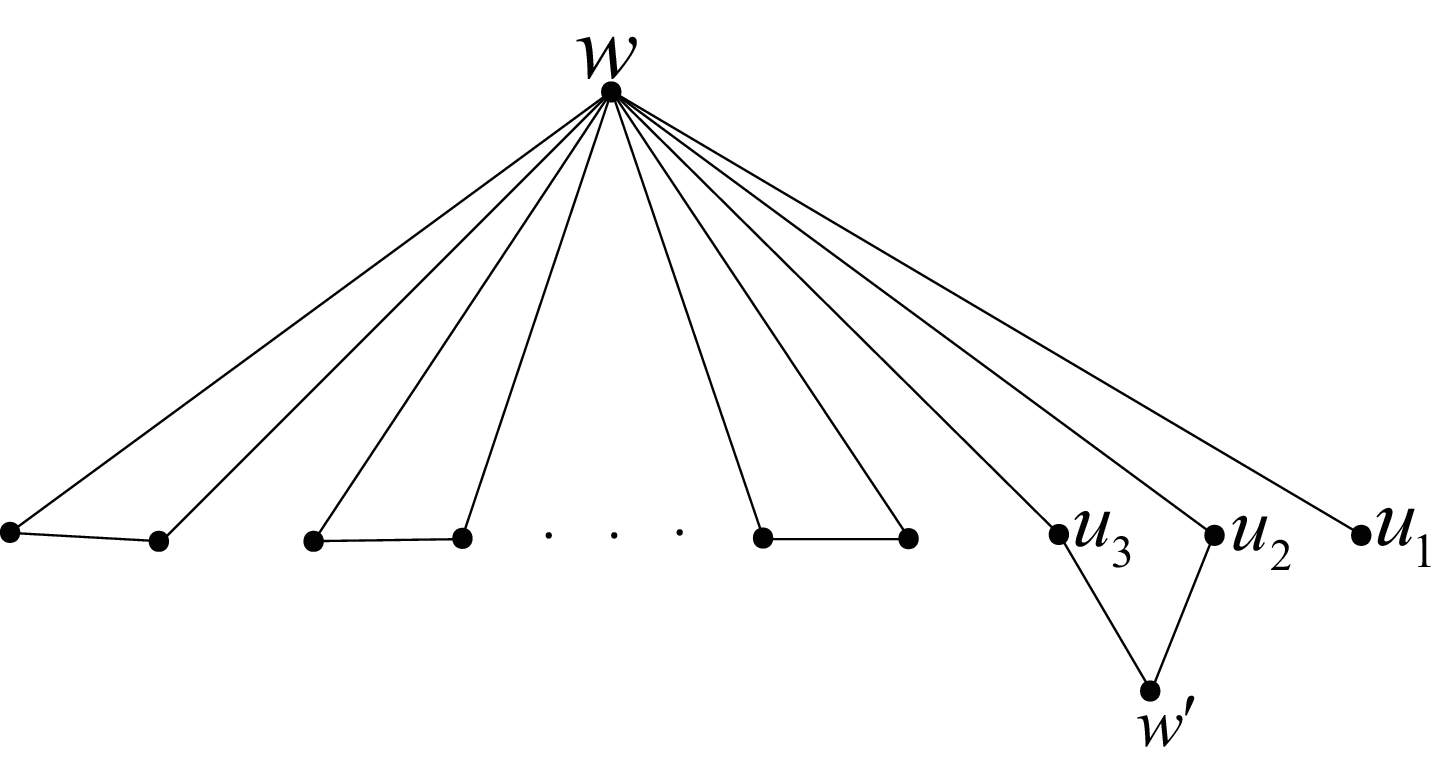}
   }
   \quad\quad
   \subfigure[$L_{2}$]{
       \includegraphics[scale=0.25]{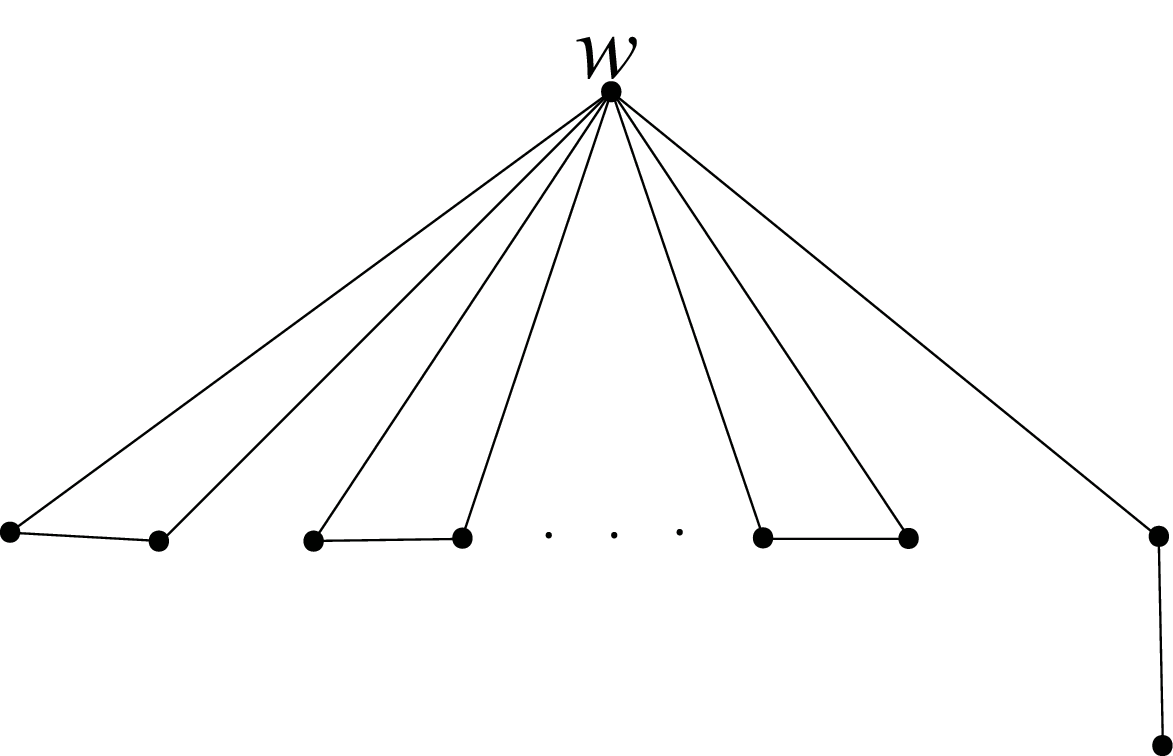}
       }
       \caption{Graphs $L_{1}$ and $L_{2}$ of Subcase 1.2.}\label{1.}
\end{centering}
\end{figure}

{\bf Case 2.} If $|V(G-w)|=2k+1$, i.e., $G\cong K_{1}\vee (G-w)$, then we consider theses cases as follows.

{\bf Subcase 2.1.} If $\delta(G)\geq2$, then $$2m=6k+4=\sum_{v\in V(G)}d(v)\geq 2k+1+2\times (2k+1)=6k+3.$$ Thus, there exists at most a vertex with degree 2 in $G-w$ and hence the degree sequence of $G-w$ is $(2, \underbrace{1,1,\cdots,1}_{2k},)$. Thus we have $G\cong K_{1}\vee ((k-1)P_{2}\cup S_{3})$.

{\bf Subcase 2.2.} If there exists a vertex with degree 1, then $$2m=6k+4=\sum_{v\in V(G)}d(v)\geq 2k+1+2\times 2k+1=6k+2.$$ Thus, there exists at most a vertex with degree 3 in $G-w$ and hence the degree sequence of $G-w$ is $(3,\underbrace{1,1,\cdots,1}_{2k-1},0)$ or $(2,2, \underbrace{1,1,\cdots,1}_{2k-2},0)$. Thus, $G\in \{ K_{1}\vee ((k-2)P_{2}\cup S_{4}\cup P_{1}), K_{1}\vee ((k-2)P_{2}\cup P_{4}\cup P_{1}), K_{1}\vee ((k-3)P_{2}\cup 2S_{3}\cup P_{1})\}$.

\noindent\begin{lemma}\label{le:ch-2.6.}{\rm(}$\cite{Wsj}${\rm)} For $k\geq2$, we have $q(K_{1}\vee ((k-1)P_{2}\cup S_{3})>2k+2+\frac{1}{k}$.
\end{lemma}

\noindent\begin{lemma}\label{le:ch-2.7.} For $k\geq2$, we have that $q(K_{1}\vee ((k-2)P_{2}\cup S_{3}\cup P_{1}))>q(K_{1}\vee kP_{2})$, $q(K_{1}\vee ((k-2)P_{2}\cup S_{4}\cup P_{1}))>q(K_{1}\vee ((k-3)P_{2}\cup 2S_{3}\cup P_{1}))>q(K_{1}\vee ((k-1)P_{2}\cup S_{3})>q(L_{2})>q(L_{1})$ and $q(K_{1}\vee ((k-2)P_{2}\cup S_{4}\cup P_{1}))>q(K_{1}\vee ((k-2)P_{2}\cup P_{4}\cup P_{1}))$.
\end{lemma}

\noindent{\bf Proof.} For convenience, Let $G_{1}=K_{1}\vee kP_{2}$ and $G_{1}^{\prime}=K_{1}\vee ((k-2)P_{2}\cup S_{3}\cup P_{1})$. Let $V(G_{1})=\{w, u_{1}, v_{1},u_{2}, v_{2},\cdots, u_{k}, v_{k}\}$ and $X$ be the perron vector $Q(G_{1})$. It is obviously that $x_{u_{i}}=x_{v_{j}}=x_{u_{1}}$ for $i,j\in \{1,2, \cdots,k\}$. Since $G_{1}^{\prime}=G_{1}-u_{1}v_{1}+u_{2} v_{1}$, we obtain $q(G_{1}^{\prime})-q(G_{1})\geq X^{T}(Q(G_{1}^{\prime})-Q(G_{1}))X=(x_{u_{2}}+x_{v_{1}})^{2}-(x_{u_{1}}+x_{v_{1}})^{2}=0$. Suppose that $q(G_{1}^{\prime})=q(G_{1})$, we have  $X$ be the perron vector $Q(G_{1}^{\prime})$ and $q(G_{1}^{\prime})x_{u_{2}}=3x_{u_{2}}+x_{v_{1}}+x_{w}+x_{v_{2}}$ and $q(G_{1})x_{u_{2}}=2x_{u_{2}}+x_{w}+x_{v_{2}}$. Moreover, $x_{u_{2}}+x_{v_{1}}=0$, which contradicts $x_{u_{2}}=\frac{x_{w}}{q(G_{1})-3}>0$. Thus, $q(G_{1}^{\prime})>q(G_{1})$.

Note that $L_{2}=L_{1}-u_{2}w^{\prime}+u_{1}u_{2}$. Let $X$ be the perron vector $Q(L_{1})$. Since $(q(L_{1})-1)x_{u_{1}}=x_{w}$,  $(q(L_{1})-2)x_{w^{\prime}}=2x_{u_{2}}$ and $(q(L_{1})-2)x_{u_{2}}=x_{w}+x_{w^{\prime}}$, we have $x_{w^{\prime}}=\frac{2x_{w}}{(q(L_{1})-2)^{2}-2}<\frac{x_{w}}{(q(L_{1})-1)}$ for $q^{2}-6q+4>0$. Thus,
$q(L_{2})-q(L_{1})\geq X^{T}(Q(L_{2})-Q(L_{1}))X=(x_{u_{1}}+x_{u_{2}})^{2}-(x_{u_{2}}+x_{w^{\prime}})^{2}>0$ , that is $q(L_{2})>q(L_{1})$. The quotient matrix $B_{\Pi}$ of the $Q$-matrix $Q(L_{2})$ is
$$
B_{\Pi}=
\left(\begin{array}{cccc}
x-2\,k-1 & -2\,k & -1 & 0\\
-1 & x-3 & 0 & 0\\
-1 & 0 & x-2 & -1\\
0 & 0 & -1 & x-1
\end{array}\right)
.$$ Then the characteristic polynomial of $B_{\Pi}$ is as follows: $$f(x)=x^{4}-(2k+7)x^{3}+(10k+15)x^{2}-(14k+9)x+4k.$$ By lemma \ref{le:ch-2.6.}, we have $q(K_{1}\vee ((k-1)P_{2}\cup S_{3})>2k+2+\frac{1}{k}$. By Matlab, we have $f(2k+2+\frac{1}{k})=\frac{2k^{5} +4k{^4} +3k^{3} +3k^{2} +k+1}{k^{4}}>0$. Thus, $q(K_{1}\vee ((k-1)P_{2}\cup S_{3})>q(L_{2})$.

Let $G_{2}=K_{1}\vee ((k-1)P_{2}\cup S_{3})$, $G_{2}^{\prime}=K_{1}\vee ((k-3)P_{2}\cup 2S_{3}\cup P_{1})$ and $V(G_{2})=\{w, u_{1}, v_{1},u_{2}, v_{2},\cdots, u_{2k-2}, v_{2k-2}, u_{1}^{\prime}, v_{1}^{\prime},v_{2}^{\prime},\}$. Let $X$ be the perron vector $Q(G_{2})$. Note that $G_{2}^{\prime}=G_{2}-u_{1}v_{1}+u_{1}u_{2}$. Since $x_{u_{i}}=x_{v_{j}}=x_{u_{1}}$ for $i,j\in \{1,2, \cdots,2k-2\}$, we have $q(G_{2}^{\prime})-q(G_{2})\geq X^{T}(Q(G_{2}^{\prime})-Q(G_{2}))X=(x_{u_{1}}+x_{u_{2}})^{2}-(x_{u_{1}}+x_{v_{1}})^{2}=0$. Suppose that $q(G_{2}^{\prime})=q(G_{2})$, we have $X$ be the perron vector $Q(G_{2}^{\prime})$ and $q(G_{2}^{\prime}-2)x_{u_{2}}=x_{u_{1}}+x_{u_{2}}+x_{v_{2}}+x_{w}=q(G_{2}-2)x_{u_{2}}=x_{v_{2}}+x_{w}$, and hence $x_{u_{2}}+x_{u_{1}}=0$, which contradicts $x_{u_{2}}=\frac{x_{w}}{q(G_{2})-3}>0$. Thus, $q(G_{2}^{\prime})>q(G_{2})$, that is $q(K_{1}\vee ((k-3)P_{2}\cup 2S_{3}\cup P_{1}))>q(K_{1}\vee ((k-1)P_{2}\cup S_{3})$.

By corollary \ref{le:ch-2.10.}, since $(k-2)P_{2}\cup S_{4}, (k-3)P_{2}\cup 2S_{3}, (k-2)P_{2}\cup P_{4}\in \mathcal{F}_{k+1,k-1}$, we have $q(K_{1}\vee ((k-2)P_{2}\cup S_{4}\cup P_{1}))>q(K_{1}\vee ((k-3)P_{2}\cup 2S_{3}\cup P_{1}))$ and $q(K_{1}\vee ((k-2)P_{2}\cup S_{4}\cup P_{1}))>K_{1}\vee ((k-2)P_{2}\cup P_{4}\cup P_{1})$. $\qedsymbol$

\noindent\begin{lemma}\label{le:ch-2.8.} For $k\geq2$, we have that $q(K_{1}\vee ((k-2)P_{2}\cup S_{3}\cup P_{1}))>2k+1+\frac{2}{2k-1}$, $q(K_{1}\vee ((k-2)P_{2}\cup P_{1})>2k+2$ and $q(K_{1}\vee ((k-2)P_{2}\cup S_{4}\cup P_{1}))>2k+2+\frac{1}{k}$.
\end{lemma}

\noindent{\bf Proof.} The quotient matrix $B_{\Pi}$ of the $Q$-matrix $Q(K_{1}\vee ((k-2)P_{2}\cup S_{3}\cup P_{1}))$ is
$$
\left(\begin{array}{ccccc}
x-2\,k & 4-2\,k & -2 & -1 & -1\\
-1 & x-3 & 0 & 0 & 0\\
-1 & 0 & x-2 & -1 & 0\\
-1 & 0 & -2 & x-3 & 0\\
-1 & 0 & 0 & 0 & x-1
\end{array}\right)
.$$ Then the characteristic polynomial of $B_{\Pi}$ is as follows: $f_{1}(x,k)=x^{5}-(2k+9)x^{4}+(16k+27)x^{3}-(42k+3)x^{2}+(44k+4)x-16k+8$ and $f_{1}(2k+1+\frac{2}{2k-1},k)=-\frac{16}{(2k-1)^{5}}(2\,k^2 -k+1)(8\,k^4 -12\,k^3 +2\,k^2 +5\,k-4)<0$ for $k\geq2$. Thus, $q(K_{1}\vee ((k-2)P_{2}\cup S_{3}\cup P_{1}))>2k+1+\frac{2}{2k-1}$. The quotient matrix $B_{\Pi}$ of the $Q$-matrix $Q(K_{1}\vee ((k-2)P_{2}\cup P_{1}))$ is
$$
B_{\Pi}=\left(\begin{array}{ccc}
x-2k-1 & -2k & -1\\
-1 & x-3 & 0\\
-1 & 0 & x-1
\end{array}\right)
.$$ Then the characteristic polynomial of $B_{\Pi}$ is as follows: $g(x,k)=(x-2)(x^{2}-(2k+3)x+2k)$. Thus, $q(K_{1}\vee (kP_{2}\cup P_{1}))=\frac{2k+3+\sqrt{4k^{2}+4k+9}}{2}>2k+2$. By Lemmas \ref{le:ch-2.7.} and \ref{le:ch-2.6.}, we have $q(K_{1}\vee ((k-2)P_{2}\cup S_{4}\cup P_{1}))>q(K_{1}\vee ((k-1)P_{2}\cup S_{3})>2k+2+\frac{1}{k}$.

Let $G$ be the 2-leaves free graph with size $m\geq17$ with maximum $Q$-index. By Lemma \ref{le:ch-2.8.}, we have
\begin{equation}\label{eq:ch-1}
q(G)\geq \left\{
\begin{aligned}
q(K_{1}\vee ((k-2)P_{2}\cup S_{3}\cup P_{1}))&>2k+1+\frac{2}{2k-1}&    &if \ m=3k;\\
q(K_{1}\vee ((k-2)P_{2}\cup P_{1}))&>2k+2 &                            &if \ m=3k+1;\\
q(K_{1}\vee ((k-2)P_{2}\cup S_{4}\cup P_{1}))&>2k+2+\frac{1}{k}&       &if \ m=3k+2.
\end{aligned}
\right.
\end{equation}
Let $w$ be a vertex of $G$ such that
$$max_{v\in V(G)}\{d(v)+m(v)\}=d(w)+m(w).$$
Note that we only need to consider the case $d(w)>2$. Otherwise, $G$ is 2-regular and thus $q(G)=4$, contradicts the
inequation \eqref{eq:ch-1}.
Recall that $$m(w)\leq min\{\frac{2m}{d(w)}-1, \frac{m}{d(w)}+\frac{d(w)-1}{2}\}.$$ By Lemma \ref{le:ch-2.3.}, we have that
$$q(G)\leq min\left\{\{d(w)+\frac{2m}{d(w)}-1,\frac{3d(w)-1}{2}+\frac{m}{d(w)}\} \right\}.$$
Let $h(x)=x+\frac{2m}{x}-1$, $z(x)=\frac{3x-1}{2}+\frac{m}{x}$. Then by the monotonicity of $h(x)$ and $z(x)$, we have that for $x\in[3, \Delta(G)]$, $$h(x)\leq max\{h(3), h(\Delta(G))\},\quad\quad z(x)\leq max\{z(3),z(\Delta(G))\}.$$
Next we show that $d(w)\notin[3, \lfloor\frac{2m+1}{3}\rfloor-2]$.

If $3\leq d(w)\leq 4$, we have that $$q(G)\leq z(3)=4+\frac{m}{3},$$ contradicts the inequation \eqref{eq:ch-1}.

If $5\leq d(w)\leq \lfloor\frac{2m+1}{3}\rfloor-2$, we have that $$q(G)\leq h(\lfloor\frac{2m+1}{3}\rfloor-2)=\lfloor\frac{2m+1}{3}\rfloor-3+\frac{2m}{\lfloor\frac{2m+1}{3}\rfloor-2},$$ and
$$
\lfloor\frac{2m+1}{3}\rfloor-3+\frac{2m}{\lfloor\frac{2m+1}{3}\rfloor-2}=
\left\{
\begin{aligned}
2k+\frac{3}{k-1}&<2k+1+\frac{2}{2k-1}&      &if\ m=3k, k\geq3;\\
2k+1+\frac{5}{2k-1}&<2k+2&                  &if\ m=3k+1, k\geq3;\\
2k+1+\frac{7}{2k-1}&<2k+2+\frac{1}{k}&      &if\ m=3k+2, k\geq3,
\end{aligned}
\right.$$ contracts the inequation \eqref{eq:ch-1}.

Thus, we only consider the case that $d(w)\in \{\lfloor\frac{2m+1}{3}\rfloor-1,\lfloor\frac{2m+1}{3}\rfloor\}$. Furthermore, if $d(w)=\lfloor\frac{2m+1}{3}\rfloor$, then by Lemma \ref{le:ch-2.5.}, we can get our conclusion immediately. So in what follows, we assume that $d(w)=\lfloor\frac{2m+1}{3}\rfloor-1$.

{\bf Case 1.} If $m=3k$, then $d(w)=\Delta(G)=2k-1$.

{\bf Subcase 1.1.} $G\cong K_{1}\vee H$ with $V(H)=2k-1$ and $|E(H)|=k+1$.

{\bf Subcase 1.1.1.} $H$ is a forest.

If $\delta(H)\geq1$, then by Lemma \ref{le:ch-2.9.}, we have $H\cong (k-3)P_{2}\cup S_{5}$ and $G\cong K_{1}\vee ((k-3)P_{2}\cup S_{5})$. We claim that $q(G)<q(K_{1}\vee ((k-2)P_{2}\cup S_{3}\cup P_{1}))$ to get a contradiction. Let $V(S_{5})=\{u,u_{1},u_{2},u_{3},u_{4}\}$ and $u$ be the center vertex of $S_{5}$. Let $X$ be the perron vector of $Q(G)$, we have $x_{u_{1}}=\cdots=x_{u_{4}}$ and $$(q(G)-5)x_{u}=x_{w}+4x_{u_{1}}, (q(G)-2)x_{u_{1}}=x_{u}+x_{w}.$$
Thus, $$x_{u}=\frac{q(G)+2}{(q(G)-1)(q(G)-6)}x_{w},\quad\quad x_{u_{1}}=\frac{q(G)-4}{(q(G)-1)(q(G)-6)}x_{w}.$$ Let $G_{3}=G\cup \{v\}$  and $G_{3}^{\prime}=G_{3}-uu_{1}-uu_{2}+u_{1}u_{2}+wv$. It is obvious that $G_{3}^{\prime}=K_{1}\vee ((k-2)P_{2}\cup S_{3}\cup P_{1})$ and $X^{\prime}=(X^{T},0)^{T}$ is an eigenvector of $Q(G_{3})$ corresponding to $q(G_{4})=q(G)$. Thus, we have
$$
\begin{aligned}
q(G_{3}^{\prime})-q(G_{3})&\geq {X^{\prime}}^{T}(Q(G_{3}^{\prime})-Q(G_{3}))X^{\prime}\\&
=(x_{w}+x_{v})^{2}+(x_{u_{1}}+x_{u_{2}})^{2}-(x_{u}+x_{u_{1}})^{2}-(x_{u}+x_{u_{1}})^{2}\\&
=\left(1-\frac{2(q(G)+2)^{2}}{(q(G)-1)^{2}(q(G)-6)^{2}}-\frac{4(q(G)+2)(q(G)-4)}{(q(G)-1)^{2}(q(G)-6)^{2}}\right)x_{w}^{2}\\&
>0,
\end{aligned}
$$
for $k\geq4$, which implies $q(G)=q(G_{3})<q(G_{3}^{\prime})$, a contradiction.

If there exists a vertex $w^{\prime}$ with $d(w^{\prime})=1$, then $$2m=6k\geq 2k-1+1+2\times(2k-1-1)=6k-4,$$  Furthermore, there exist at most a vertex with degree 6 in $G$. We consider these cases in the following.

If $$2m=6k=2k-1+1+2x+3(2k-1-1-x),$$  then $x=2k-6$ and $2k-1-1-x=4$ and hence the degree sequence of $G-w$ is $\{2,2,2,2,\underbrace{1,1,\cdots, 1}_{2k-6},0\}$. Thus, $G\in \{K_{1}\vee ((k-4)P_{2}\cup P_{6}\cup P_{1}), (k-5)P_{2}\cup P_{5}\cup S_{3}\cup P_{1}),(k-5)P_{2}\cup 2P_{4}\cup P_{1}),(k-7)P_{2}\cup 4S_{3}\cup P_{1})\}$.

If $$2m=6k=2k-1+1+2x+4(2k-1-1-x),$$  then $x=2k-4$ and $2k-1-1-x=2$ and hence the degree sequence of $G-w$ is $\{3,3,3,\underbrace{1,1,\cdots, 1}_{2k-4},0\}$. Thus, $G\in \{K_{1}\vee ((k-4)P_{2}\cup D_{2,2}\cup P_{1}), (k-5)P_{2}\cup 2S_{4}\cup P_{1})\}$.

If $$2m=6k=2k-1+1+2x+5(2k-1-1-x),$$ then $x=\frac{6k-10}{3}$, a contradiction.

If $$2m=6k=2k-1+1+2x+6(2k-1-1-x),$$ then $x=2k-3$ and $2k-1-1-x=1$ and hence the degree sequence of $G-w$ is $\{5,\underbrace{1,1,\cdots, 1}_{2k-3},0\}$. Thus, $G\in \{K_{1}\vee ((k-4)P_{2}\cup S_{6}\cup P_{1})\}$.

If $$2m=6k=2k-1+1+2x+7(2k-1-1-x),$$ then $x=\frac{10k-13}{5}$, a contradiction.

If $$2m=6k=2k-1+1+3x+4(2k-1-1-x),$$ then $x=4k-8<2k-2$, a contradiction.

Note that $\{K_{1}\vee ((k-4)P_{2}\cup P_{6}\cup P_{1}), (k-5)P_{2}\cup P_{5}\cup S_{3}\cup P_{1}),(k-5)P_{2}\cup 2P_{4}\cup P_{1}),(k-7)P_{2}\cup 4S_{3}\cup P_{1}), K_{1}\vee ((k-4)P_{2}\cup D_{2,2}\cup P_{1}), (k-5)P_{2}\cup 2S_{4}\cup P_{1}), K_{1}\vee ((k-4)P_{2}\cup S_{6}\cup P_{1})\}\in K_{1}\vee (F\cup P_{1})$, where $F\in\mathcal{F}_{k+1,k-3}$. By Corollary \ref{le:ch-2.10.}, we have $G\cong K_{1}\vee ((k-4)P_{2}\cup S_{6}\cup P_{1})$. If $G\cong K_{1}\vee ((k-4)P_{2}\cup S_{6}\cup P_{1})$, then let $V(S_{6})=\{u, u_{1},u_{2},u_{3},u_{4},u_{5}\}$ and $u$ be the center vertex of $S_{6}$. Let $X$ be the perron vector of $Q(G)$, we have $x_{u_{1}}=\cdots=x_{u_{5}}$ and $$(q(G)-6)x_{u}=x_{w}+5x_{u_{1}}, \quad\quad (q(G)-2)x_{u_{1}}=x_{u}+x_{w}.$$ Thus, $$x_{u}=\frac{q(G)+3}{(q(G)-1)(q(G)-7)}x_{w},\quad\quad x_{u_{1}}=\frac{q(G)-5}{(q(G)-1)(q(G)-7)}x_{w}.$$ Let $G_{4}=G\cup \{v\}$  and $G_{4}^{\prime}=G_{6}-uu_{1}-uu_{2}-uu_{3}+u_{1}u_{2}+u_{3}v+wv$. It is obvious that $G_{4}^{\prime}=K_{1}\vee ((k-2)P_{2}\cup S_{3}\cup P_{1})$ and $X^{\prime}=(X^{T},0)^{T}$ is an eigenvector of $Q(G_{4})$ corresponding to $q(G_{4})=q(G)$. Thus, we have
$$
\begin{aligned}
q(G_{4}^{\prime})-q(G_{4})&\geq {X^{\prime}}^{T}(Q(G_{4}^{\prime})-Q(G_{4}))X^{\prime}\\
&=(x_{w}+x_{v})^{2}+(x_{u_{1}}+x_{u_{2}})^{2}+(x_{u_{3}}+x_{v})^{2}-(x_{u}+x_{u_{1}})^{2}-(x_{u}+x_{u_{1}})^{2}\\
&-(x_{u}+x_{u_{3}})^{2}\\
&=(1-\frac{2(q(G)-5)^{2}}{(q(G)-1)^{2}(q(G)-7)^{2}}-\frac{3(q(G)+3)^{2}}{(q(G)-1)^{2}(q(G)-7)^{2}}\\
&-\frac{6(q(G)+3)(q(G)-5)}{(q(G)-1)^{2}(q(G)-7)^{2}})x_{w}^{2}>0,
\end{aligned}
$$
for $k\geq5$, which implies $q(G)=q(G_{4})<q(K_{1}\vee ((k-2)P_{2}\cup S_{3}\cup P_{1}))$, a contradiction.

{\bf Subcase 1.1.2.} $H$ contains a cycle $u_{1}u_{2}\cdots u_{s}u_{1}$.

If $\delta(G)\geq2$, then $$2k+2=2|E(H)|\geq \sum_{i=1}^{s}d_{H}(u_{i})+(2k-1-s)\times 1=2k+s-1.$$ Thus, $s=3,d_{H}(u_{i})=2$ for $i=1,2,3$, $d_{H}(v)=1$ for $v\in V(H)\setminus \{u_{1},u_{2},u_{3}\}$. Hence we have $H=(k-2)P_{2}\cup C_{3}$ and $G\cong K_{1}\vee ((k-2)P_{2}\cup C_{3})$. We claim that $q(G)<q(K_{1}\vee ((k-2)P_{2}\cup S_{3}\cup P_{1}))$ to get a contradiction. Let $X$ be the perron vector of $Q(G)$, we have $x_{u_{1}}=x_{u_{2}}=x_{u_{3}}$ and
$$(q(G)-3)x_{u_{1}}=x_{w}+2x_{u_{1}}.$$
Thus, $$x_{u_{1}}=\frac{1}{q(G)-5}x_{w}.$$ Let $G_{5}=G\cup \{v\}$ and $G_{5}^{\prime}=G_{5}-u_{1}u_{2}+wv$. It is obvious that $G_{5}^{\prime}=K_{1}\vee ((k-2)P_{2}\cup S_{3}\cup P_{1})$ and $X^{\prime}=(X^{T},0)^{T}$ is an eigenvector of $Q(G_{5})$ corresponding to $q(G_{5})=q(G)$. Thus, we have
$$
\begin{aligned}
q(G_{5}^{\prime})-q(G_{5})&\geq {X^{\prime}}^{T}(Q(G_{5}^{\prime})-Q(G_{5}))X^{\prime}\\
&=(x_{w}+x_{v})^{2}-(x_{u_{1}}+x_{u_{2}})^{2}\\
&=\left(1-\frac{4}{(q(G)-5)^{2}}\right)x_{w}^{2}>0,
\end{aligned}
$$
for $x>2k+1+\frac{1}{2k-1}>7,k\geq3$, which implies $q(G)=q(G_{5})<q(K_{1}\vee ((k-2)P_{2}\cup S_{3}\cup P_{1}))$, a contradiction. If there exists a vertex $w^{\prime}$ with $d(w^{\prime})=1$, then $$2k+2=2|E(H)|\geq \sum_{i=1}^{s}d_{H}(u_{i})+(2k-2-s)\times 1=2k+s-2.$$ Thus, $s\leq4$. Suppose that $s=4$, we have $d_{H}(u_{i})=2$ for $i=1,2,3,4$, $d_{H}(v)=1$ for $v\in V(H)\setminus \{u_{1},u_{2},u_{3},w^{\prime}\}$. Hence we have $H\cong(k-3)P_{2}\cup C_{4}\cup P_{1}$ and $G\cong K_{1}\vee ((k-3)P_{2}\cup C_{4}\cup P_{1})$. We claim that $q(G)<q(K_{1}\vee ((k-2)P_{2}\cup S_{3}\cup P_{1}))$ to get a contradiction. Let $X$ be the perron vector of $Q(G)$, we have $x_{u_{1}}=x_{u_{2}}=x_{u_{3}}=x_{u_{4}}$ and
$$(q(G)-3)x_{u_{1}}=x_{w}+2x_{u_{1}}.$$
Thus, $$x_{u_{1}}=\frac{1}{q(G)-5}x_{w}.$$ Let $G_{6}=G\cup \{v\}$ and $G_{6}^{\prime}=G_{6}-u_{1}u_{2}-u_{1}u_{3}+u_{1}v+wv$. It is obvious that $G_{6}^{\prime}=K_{1}\vee ((k-2)P_{2}\cup S_{3}\cup P_{1})$ and $X^{\prime}=(X^{T},0)^{T}$ is an eigenvector of $Q(G_{6})$ corresponding to $q(G_{6})=q(G)$. Thus, we have
$$
\begin{aligned}
q(G_{6}^{\prime})-q(G_{6})&\geq {X^{\prime}}^{T}(Q(G_{6}^{\prime})-Q(G_{6}))X^{\prime}\\
&=(x_{w}+x_{v})^{2}+(x_{u_{1}}+x_{v})^{2}-(x_{u_{1}}+x_{u_{2}})^{2}-(x_{u_{1}}+x_{u_{3}})^{2}\\
&={x_{w}}^{2}-7x_{u_{1}}^{2}\\
&=\left(1-\frac{7}{(q(G)-5)^{2}}\right)x_{w}^{2}>0,
\end{aligned}
$$
for $x>2k+1+\frac{1}{2k-1}>9,k\geq4$, which implies $q(G)=q(G_{6})<q(K_{1}\vee ((k-2)P_{2}\cup S_{3}\cup P_{1}))$, a contradiction. Suppose that $s=3$, the degree sequence of $G-w$ is $\{2,2,2,2,\underbrace{1,1,\cdots,1}_{k-4},0\}$. Thus, $H\cong (k-4)P_{2}\cup S_{3}\cup C_{3}\cup P_{1}$ and $G\cong K_{1}\vee ((k-4)P_{2}\cup S_{3}\cup C_{3}\cup P_{1})$. We claim that $q(G)<q(K_{1}\vee ((k-2)P_{2}\cup S_{3}\cup P_{1}))$ to get a contradiction. Let $X$ be the perron vector of $Q(G)$, we have $x_{u_{1}}=x_{u_{2}}=x_{u_{3}}$ and
$$(q(G)-3)x_{u_{1}}=x_{w}+2x_{u_{1}}.$$
Thus, $$x_{u_{1}}=\frac{1}{q(G)-5}x_{w}.$$ Let $G_{7}=G\cup \{v\}$ and $G_{7}^{\prime}=G_{7}-u_{1}u_{2}-u_{1}u_{3}+u_{1}v+wv$. It is obvious that $G_{7}^{\prime}=K_{1}\vee ((k-2)P_{2}\cup S_{3}\cup P_{1})$ and $X^{\prime}=(X^{T},0)^{T}$ is an eigenvector of $Q(G_{7})$ corresponding to $q(G_{7})=q(G)$. Thus, we have
$$
\begin{aligned}
q(G_{7}^{\prime})-q(G_{7})&\geq {X^{\prime}}^{T}(Q(G_{7}^{\prime})-Q(G_{7}))X^{\prime}\\
&=(x_{w}+x_{v})^{2}+(x_{u_{1}}+x_{v})^{2}-(x_{u_{1}}+x_{u_{2}})^{2}-(x_{u_{1}}+x_{u_{3}})^{2}\\
&={x_{w}}^{2}-7x_{u_{1}}^{2}\\
&=\left(1-\frac{7}{(q(G)-5)^{2}}\right)x_{w}^{2}>0,
\end{aligned}
$$
for $x>2k+1+\frac{1}{2k-1}>9,k\geq4$, which implies $q(G)=q(G_{7})<q(K_{1}\vee ((k-2)P_{2}\cup S_{3}\cup P_{1}))$, a contradiction.

{\bf Subcase 1.2.} There exists $v\in V(G)$ such that $v\notin N_{G}(w)$.

If $\delta(G)\geq2$, then $$\sum_{u\in N_{G}(w)}d_{G}(u)\leq 2\times 3k-(2k-1)-d_{G}(v)\leq 4k-1.$$ Then we have $$q(G)\leq d_{G}(w)+m(w)\leq 2k-1+\frac{4k-1}{2k-1}=2k+1+\frac{1}{2k-1},$$ contradicts the inequation \eqref{eq:ch-1}. If there exists a vertex $w^{\prime}$ with $d(w^{\prime})=1$, then $$\sum_{u\in N_{G}(w)}d_{G}(u)\leq 2\times 3k-(2k-1)-d_{G}(v)\leq 4k.$$ Then we have $$q(G)\leq d_{G}(w)+m(w)\leq 2k-1+\frac{4k}{2k-1}=2k+1+\frac{2}{2k-1},$$ contradicts the inequation \eqref{eq:ch-1}.

{\bf Case 2.} $m=3k+1$. Then $d(w)=\Delta(G)=2k$.

{\bf Subcase 2.1.} $G=K_{1}\vee H$ with $|V(H)|=2k$ and $|E(H)|=k+1$.

{\bf Subcase 2.1.1.} $H$ is a forest.

If $\delta(G)\geq2$, then by Lemma \ref{le:ch-2.9.}, we have $H\in \mathcal{F}_{k+1,k-1}$, that is $H\cong S_{4}\cup (k-2)P_{2}$ and hence $G\cong K_{1}\vee ((k-2)P_{2}\cup S_{4})$. Let $V(S_{4})=\{u, u_{1},u_{2},u_{3}\}$ and $u$ be the center vertex of $S_{4}$. Let $X$ be the perron vector of $Q(G)$, we have $x_{u_{1}}=x_{u_{2}}=x_{u_{3}}$ and $$(q(G)-4)x_{u}=x_{w}+3x_{u_{1}}, \quad\quad (q(G)-2)x_{u_{1}}=x_{u}+x_{w}.$$
Thus, $$x_{u}=\frac{q(G)+1}{(q(G)-1)(q(G)-5)}x_{w},\quad\quad x_{u_{1}}=\frac{q(G)-3}{(q(G)-1)(q(G)-5)}x_{w}.$$ Let $G_{8}=G\cup \{v\}$ and $G_{8}^{\prime}=G_{8}-uu_{1}-uu_{2}-uu_{3}+u_{1}u_{2}+u_{3}v+wv$. It is obvious that $G_{8}^{\prime}=K_{1}\vee(kP_{2}\cup P_{1})$ and $X^{\prime}=(X^{T},0)^{T}$ is an eigenvector of $Q(G_{8})$ corresponding to $q(G_{8})=q(G)$. Thus, we have
$$
\begin{aligned}
q(G_{8}^{\prime})-q(G_{8})&\geq {X^{\prime}}^{T}(Q(G_{8}^{\prime})-Q(G_{8}))X^{\prime}\\
&=(x_{w}+x_{v})^{2}+(x_{u_{1}}+x_{u_{2}})^{2}+(x_{u_{3}}+x_{v})^{2}-(x_{u}+x_{u_{1}})^{2}\\
&-(x_{u}+x_{u_{2}})^{2}-(x_{u}+x_{u_{3}})^{2}\\
&=(1-\frac{2(q(G)-3)^{2}}{(q(G)-1)^{2}(q(G)-5)^{2}}-\frac{3(q(G)+1)^{2}}{(q(G)-1)^{2}(q(G)-5)^{2}}\\
&-\frac{6(q(G)+1)(q(G)-3)}{(q(G)-1)^{2}(q(G)-5)^{2}})x_{w}^{2}\\
&=\left(\frac{q^{4}-12q^{3}+14q^{2}-41q+58}{(q(G)-1)^{2}(q(G)-5)^{2}}\right)x_{w}^{2}>0,
\end{aligned}
$$
for $q>2k+2\geq12, k\geq5$, which implies $q(G)=q(G_{8})<q(K_{1}\vee (kP_{2}\cup P_{1}))$, a contradiction.

If there exists a vertex $w^{\prime}$ satisfying $d(w^{\prime})=1$, then $$2m=6k+2=\sum_{u\in V(G)}d_{G}(u)\geq2k+1+(2k-1)\times2=6k-1,$$ and hence there exists at most a vertex with degree 5 of $G$. We consider these cases in the following.

If $$2m=6k+2=2k+1+2x+3(2k-1-x)),$$ then $x=2k-4$ and $2k-1-x=3$ and hence the degree sequence of $G-w$ is $\{2,2,2,\underbrace{1,1,\cdots, 1}_{2k-4},0\}$. Thus, $G\in \{K_{1}\vee ((k-3)P_{2}\cup P_{5}\cup P_{1}), (k-4)P_{2}\cup P_{4}\cup S_{3}\cup P_{1}),(k-5)P_{2}\cup 3S_{3}\cup P_{1})\}$.

If $$2m=6k+2=2k+1+2x+4(2k-1-x),$$ then $x=\frac{4k-5}{2}$, a contradiction.

If $$2m=6k+2=2k+1+2x+5(2k-1-x),$$ then $x=2k-2$ and $2k-1-x=1$ and hence the degree sequence of $G-w$ is $\{4, \underbrace{1,1,\cdots, 1}_{2k-2},0\}$. Thus, $G\in \{K_{1}\vee ((k-3)P_{2}\cup S_{5}\cup P_{1})\}$.

If $$2m=6k+2=2k+1+3x+4(2k-1-x),$$ then $x=4k-5<2k-1$, a contradiction.

Note that $\{K_{1}\vee ((k-3)P_{2}\cup P_{5}\cup P_{1}), (k-4)P_{2}\cup P_{4}\cup S_{3}\cup P_{1}),(k-5)P_{2}\cup 3S_{3}\cup P_{1}),K_{1}\vee ((k-3)P_{2}\cup S_{5}\cup P_{1})\}\in K_{1}\vee(F\cup P_{1})$, where $F\in\mathcal{F}_{k+1,k-2}$. By Corollary \ref{le:ch-2.10.}, we have $G\cong K_{1}\vee ((k-3)P_{2}\cup S_{5}\cup P_{1})$. If $G\cong K_{1}\vee ((k-3)P_{2}\cup S_{5}\cup P_{1})$, then let $V(S_{5})=\{u, u_{1},u_{2},u_{3}, u_{4}\}$ and $u$ be the center vertex of $S_{5}$. Let $X$ be the perron vector of $Q(G)$, we have $x_{u_{1}}=x_{u_{2}}=x_{u_{3}}=x_{u_{4}}$ and $$(q(G)-5)x_{u}=x_{w}+4x_{u_{1}}, \quad\quad (q(G)-2)x_{u_{1}}=x_{u}+x_{w}.$$
Thus, $$x_{u}=\frac{q(G)+2}{(q(G)-1)(q(G)-6)}x_{w},\quad\quad x_{u_{1}}=\frac{q(G)-4}{(q(G)-1)(q(G)-6)}x_{w}.$$ Let $G_{9}=G\cup \{v\}$  and $G_{9}^{\prime}=G_{9}-uu_{1}-uu_{2}-uu_{3}-uu_{4}+u_{1}u_{2}+u_{3}u_{4}+uv+wv$. It is obvious that $G_{9}^{\prime}=K_{1}\vee (kP_{2}\cup P_{1})$ and $X^{\prime}=(X^{T},0)^{T}$ is an eigenvector of $Q(G_{9})$ corresponding to $q(G_{9})=q(G)$. Thus, we have
$$
\begin{aligned}
q(G_{9}^{\prime})-q(G_{9})&\geq {X^{\prime}}^{T}(Q(G_{9}^{\prime})-Q(G_{9}))X^{\prime}\\
&=(x_{w}+x_{v})^{2}+(x_{u}+x_{v})^{2}+(x_{u_{1}}+x_{u_{2}})^{2}+(x_{u_{3}}+x_{u_{4}})^{2}-(x_{u}+x_{u_{1}})^{2}\\
&-(x_{u}+x_{u_{2}})^{2}-(x_{u}+x_{u_{3}})^{2}-(x_{u}+x_{u_{4}})^{2}\\
&=x_{w}^{2}+4x_{u_{1}}^{2}-3x_{u}^{2}-8x_{u}x_{u_{1}}\\
&=(1+\frac{4(q(G)-4)^{2}}{(q(G)-1)^{2}(q(G)-6)^{2}}-\frac{3(q(G)+2)^{2}}{(q(G)-1)^{2}(q(G)-6)^{2}}\\
&-\frac{8(q(G)+2)(q(G)-4)}{(q(G)-1)^{2}(q(G)-6)^{2}})x_{w}^{2}\\
&=\frac{q^{4}-14q^{3}+54q^{2}-112q+152}{(q(G)-1)^{2}(q(G)-6)^{2}}x_{w}^{2}>0,
\end{aligned}
$$
for $q>2k+2\geq10, k\geq4$, which implies $q(G)=q(G_{9})<q(K_{1}\vee (kP_{2}\cup P_{1}))$, a contradiction.

{\bf Subcase 2.1.2.} $H$ contains a cycle $u_{1}u_{2}\cdots u_{s}u_{1}$.

If $\delta(G)\geq2$, then $$2(k+1)=2|E(H)|\geq \sum_{i=1}^{s}d_{H}(u_{i})+(2k-s)\times 1=2k+s,$$ Thus $3\leq s\leq2$, a contradiction. If there exists a vertex $w^{\prime}$ satisfying $d(w^{\prime})=1$, then $$2(k+1)=2|E(H)|\geq \sum_{i=1}^{s}d_{H}(u_{i})+(2k-1-s)\times 1=2k+s-1,$$ and hence $s\leq3$, that is $s=3$. Furthermore, $H\cong (k-2)P_{2}\cup C_{3}\cup P_{1}$ and $G\cong K_{1}\vee ((k-2)P_{2}\cup C_{3}\cup P_{1})$. then Let $V(C_{3})=\{u_{1},u_{2},u_{3}\}$ and $X$ be the perron vector of $Q(G)$, we have $x_{u_{1}}=x_{u_{2}}=x_{u_{3}}=$ and $$(q(G)-3)x_{u_{1}}=x_{w}+2x_{u_{1}}.$$
Thus, $$x_{u_{1}}=\frac{1}{(q(G)-5)}x_{w}.$$ Let $G_{10}=G\cup \{v\}$ and $G_{10}^{\prime}=G_{10}-u_{1}u_{2}-u_{1}u_{3}+u_{1}v+wv$. It is obvious that $G_{10}^{\prime}=K_{1}\vee (kP_{2}\cup P_{1})$ and $X^{\prime}=(X^{T},0)^{T}$ is an eigenvector of $Q(G_{10})$ corresponding to $q(G_{10})=q(G)$. Thus, we have
$$
\begin{aligned}
q(G_{10}^{\prime})-q(G_{10})&\geq {X^{\prime}}^{T}(Q(G_{10}^{\prime})-Q(G_{10}))X^{\prime}\\
&=(x_{w}+x_{v})^{2}+(x_{u_{1}}+x_{v})^{2}-(x_{u_{1}}+x_{u_{2}})^{2}+(x_{u_{1}}+x_{u_{3}})^{2}\\
&=x_{w}^{2}-7x_{u_{1}}^{2}\\
&=(1-\frac{7}{(q(G_{10})-5))^{2}}>0,
\end{aligned}
$$
for $q>2k+2\geq8, k\geq3$, which implies $q(G)=q(G_{10})<q(K_{1}\vee (kP_{2}\cup P_{1}))$, a contradiction.

\begin{figure}[H]
\begin{centering}
 \subfigure[$L_{3}$]{
  \includegraphics[scale=0.25]{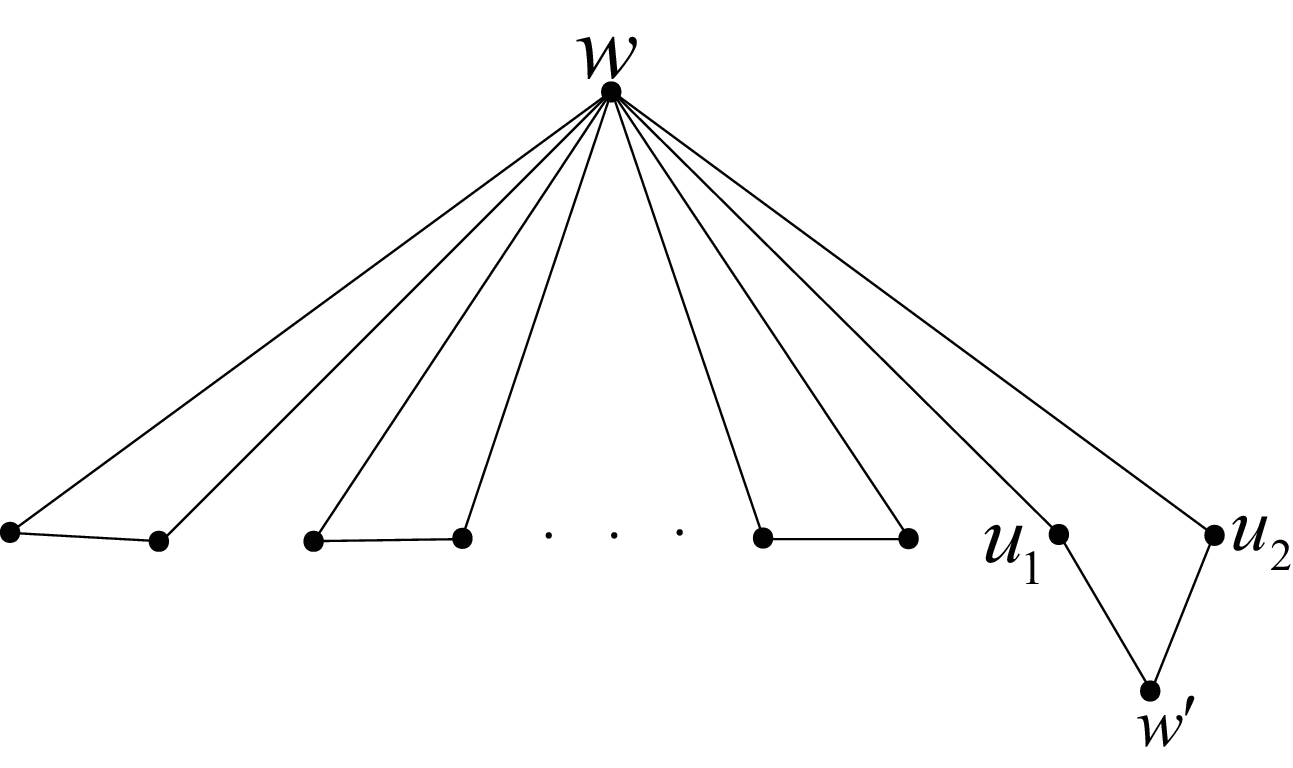}
   }
   \quad\quad
   \subfigure[$L_{4}$]{
       \includegraphics[scale=0.25]{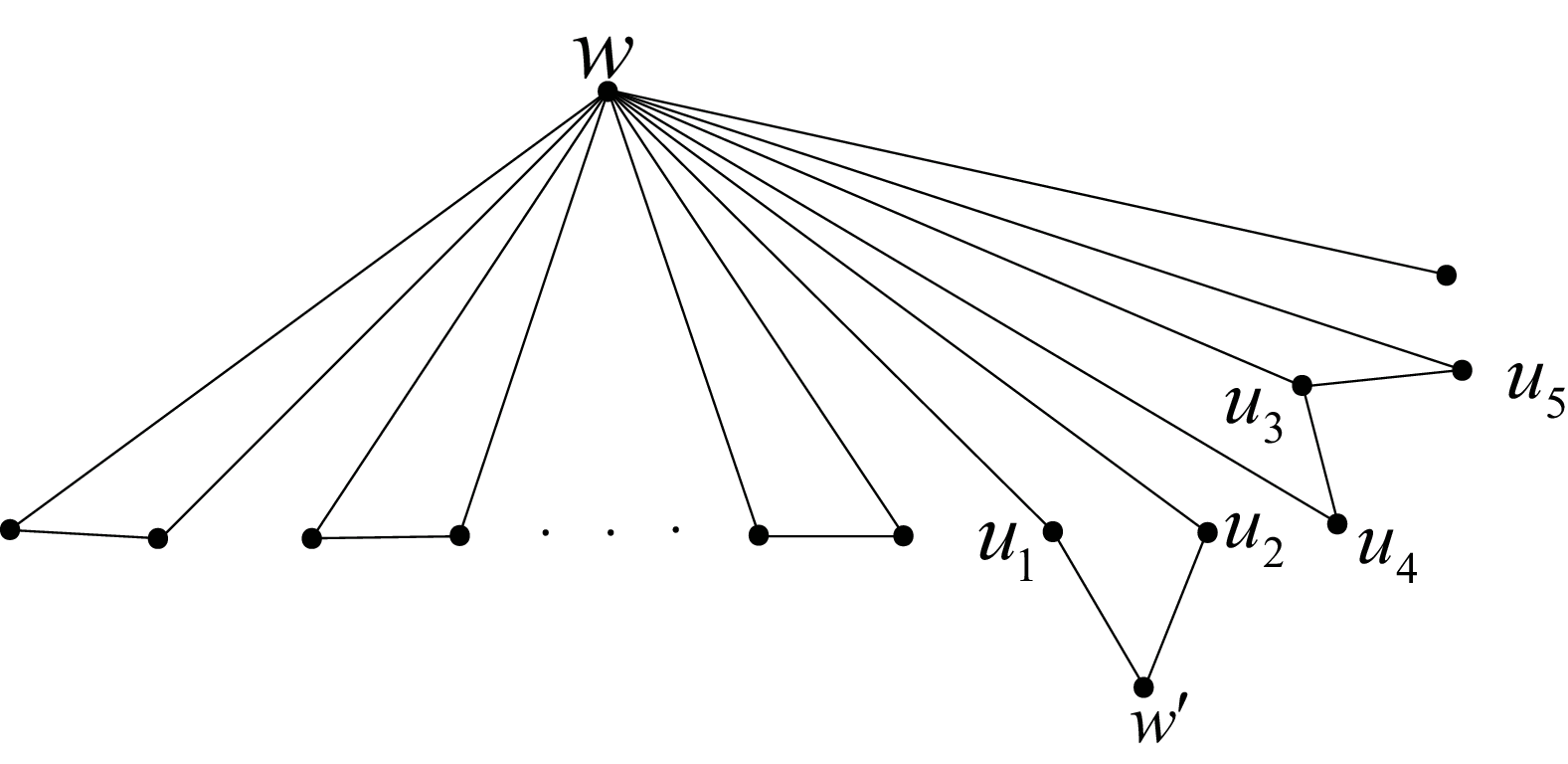}
       }
       \quad\quad
   \subfigure[$L_{5}$]{
       \includegraphics[scale=0.25]{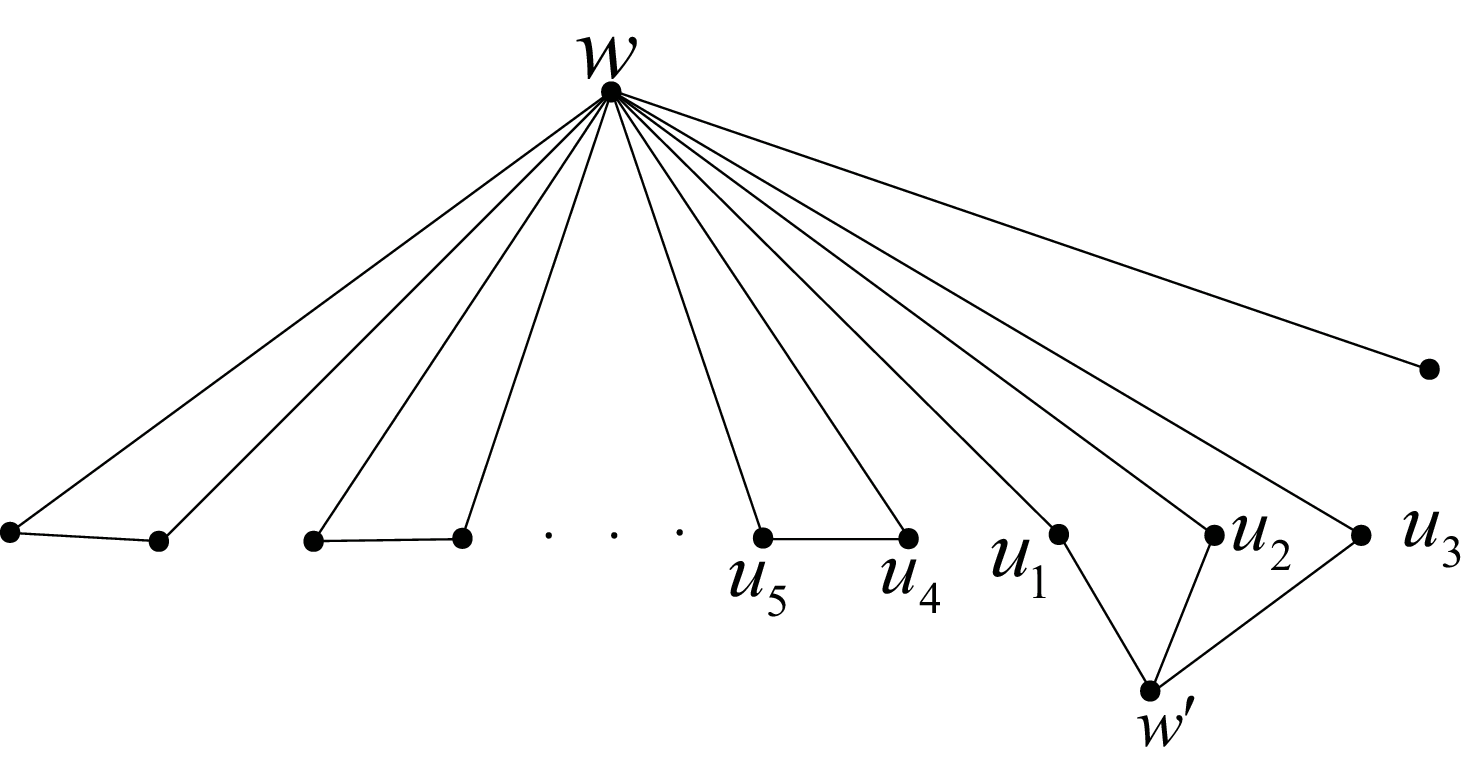}
       }
       \caption{Graphs $L_{3}$, $L_{4}$ and $L_{5}$ of Subcase 2.1.3.}\label{2.}
\end{centering}
\end{figure}

{\bf Subcase 2.1.3.} There exists $v\in V(G)$ such that $v\notin N_{G}(w)$.

If $\delta(G)\geq2$, then $$ 2(k+1)=2|E(H)|=2k+2|W|,$$ and $|W|=1$. Furthermore, $G\cong L_{3}$ (see Figure. 2). Let $X$ be the perron vector of $Q(G)$, we have $x_{u_{1}}=x_{u_{2}}$ and $$(q(G)-2)x_{u_{1}}=x_{w}+x_{w^{\prime}}, \quad (q(G)-2)x_{w^{\prime}}=2x_{u_{1}}.$$ Thus, $$x_{u_{1}}=\frac{q-2}{q(G)^{2}-4q+2}x_{w}, \quad x_{w^{\prime}}=\frac{2}{q(G)^{2}-4q+2}x_{w}.$$ Let $G^{\prime}=G-w^{\prime}u_{1}-w^{\prime}u_{2}+u_{1}u_{2}+ww^{\prime}$. It is obvious that $G^{\prime}=K_{1}\vee (kP_{2}\cup P_{1})$ and $X^{\prime}=(X^{T},0)^{T}$. Thus, we have
$$
\begin{aligned}
q(G^{\prime})-q(G)&\geq {X^{\prime}}^{T}(Q(G^{\prime})-Q(G))X^{\prime}\\
&=(x_{w}+x_{w^{\prime}})^{2}+(x_{u_{1}}+x_{u_{2}})^{2}-(x_{w^{\prime}}+x_{u_{1}})^{2}-(x_{w^{\prime}}+x_{u_{1}})^{2}\\
&=(x_{w}-x_{w^{\prime}})^{2}+2x_{u_{1}}^{2}-4x_{u_{1}}x_{w^{\prime}}\\
&=\frac{(q^{2}-4q)^{2}+2(q-2)(q-6)}{(q(G)-4q+2))^{2}}>0,
\end{aligned}
$$
for $q>2k+2\geq6, k\geq2$, which implies $q(G)<q(K_{1}\vee (kP_{2}\cup P_{1}))$, a contradiction. If there exists a vertex $w^{\prime}$ with $ d(w^{\prime})=1$, then either $$ 2(k+1)\geq 2|E(H)|=2k-1+2|W|,$$ and $|W|\leq \frac{3}{2}$  or $$2(k+1)=2|E(H)|\geq 2k+1+2(|W|-1),$$ and $|W|\leq \frac{3}{2}$. Thus $|W|=1$. Note that $2(k+1)=2|E(H)|=2k-2+2+2|W|=2k-1+3$. Thus, the degree sequence of $G-w$ is $(2,\underbrace{1,1,\cdots,1}_{2k-2},0,2)$ or $(\underbrace{1,1,\cdots,1}_{2k-1},0,3)$, and hence $H\cong L_{4}$ or $H\cong L_{5}$(see Figure. 2). Let $X$ be the perron vector of $Q(L_{5})$, we have $x_{u_{1}}=x_{u_{2}}=x_{u_{3}}, x_{u_{4}}=x_{u_{5}}$ and $$(q(G)-2)x_{u_{4}}=x_{u_{4}}+x_{w}, \quad (q(G)-3)x_{w^{\prime}}=3x_{u_{1}}, \quad(q(G)-2)x_{u_{1}}=x_{w}+x_{w^{\prime}}.$$ Thus, $x_{w^{\prime}}=\frac{3}{q^{2}-5q+3}x_{w}$ and $(q(G)-3)(x_{u_{4}}-x_{w^{\prime}})=x_{w}-3x_{u_{1}}=\frac{(q(G)-5)(q(G)^{2}-5q+3)-9}{(q(G)-2)(q(G)^{2}-5q+3)}x_{w}>0$ for
$q>2k+2\geq6, k\geq2$, which implies $q(L_{5})<q(L_{4})$. We will prove $q(L_{4})< q(K_{1}\vee (kP_{2}\cup P_{1}))$ in the following. Let $X$ be the perron vector of $Q(L_{4})$, we have $x_{u_{1}}=x_{u_{2}}, x_{u_{4}}=x_{u_{5}}$ and $$(q(G)-2)x_{u_{1}}=x_{w}+x_{w^{\prime}}, \quad (q(G)-2)x_{w^{\prime}}=2x_{u_{1}},$$

$$(q(G)-3)x_{u_{3}}=2x_{u_{4}}+x_{w}, \quad (q(G)-2)x_{u_{4}}=x_{u_{3}}+x_{w}.$$

Thus, $$x_{w^{\prime}}=\frac{2}{q(G)-2}x_{u_{1}}, \quad x_{u_{3}}=\frac{q}{(q(G)-1)(q(G)-4)}x_{w}, \quad x_{u_{4}}=\frac{q-2}{(q(G)-1)(q(G)-4)}x_{w}.$$ Let $G_{11}=K_{1}\vee (kP_{2}\cup P_{1})$. Note $G_{11}=L_{4}-u_{1}w^{\prime}-u_{2}w^{\prime}-u_{3}u_{4}+u_{1}u_{2}+u_{4}v+vw$. Let $X^{\prime}=(X^{T},0)^{T}$ is an eigenvector of $Q(G_{11})$ corresponding to $q(G_{11})=q(G)$. Thus, we have
$$
\begin{aligned}
q(G_{11})-q(L_{4})&\geq {X^{\prime}}^{T}(Q(G_{11})-Q(L_{4}))X^{\prime}\\
&=(x_{u_{4}}+x_{v})^{2}+(x_{u_{1}}+x_{u_{2}})^{2}+(x_{w}+x_{v})^{2}-(x_{w^{\prime}}+x_{u_{1}})^{2}\\
&-(x_{w^{\prime}}+x_{u_{2}})^{2}-(x_{u_{3}}+x_{u_{4}})^{2}\\
&=(2x_{u_{1}}^{2}-2x_{w^{\prime}}^{2}-4x_{u_{1}}x_{w^{\prime}})+(x_{w}^{2}-x_{u_{3}}^{2}-2x_{u_{3}}x_{u_{4}})\\
&=\frac{2q(L_{4})^{2}-16q+16}{(q(L_{4})-2)^{2}}x_{u_{1}}+\left(1-\frac{3q^{2}-4q}{(q(L_{4})-1)^{2}(q(L_{4})-4)^{2}}\right)x_{w}^{2}>0,
\end{aligned}
$$
for $q>2k+2\geq6, k\geq3$, which implies $q(L_{4})<q(K_{1}\vee (kP_{2}\cup P_{1}))$, a contradiction.

{\bf Case 3.} $m=3k+2$. Then $d(w)=\Delta(G)=2k$.

{\bf Subcase 3.1.} $G=K_{1}\vee H$ with $|V(H)|=2k$ and $|E(H)|=k+2$.

{\bf Subcase 3.1.1.} $H$ is a forest.

If $\delta(G)\geq2$, then by Lemma \ref{le:ch-2.9.}, we have $H\in \mathcal{F}_{k+2,k-2}$, that is $H\cong S_{6}\cup (k-3)P_{2}$ and hence $G\cong K_{1}\vee ((k-3)P_{2}\cup S_{6})$. If $G\cong K_{1}\vee ((k-3)P_{2}\cup S_{6})$, then let $V(S_{6})=\{u, u_{1},u_{2},u_{3}, u_{4},u_{5}\}$ and $u$ be the center vertex of $S_{6}$. Let $X$ be the perron vector of $Q(G)$, we have $x_{u_{1}}=x_{u_{i}}$ for $i\in \{2,3,4,5\}$ and $$(q(G)-6)x_{u}=x_{w}+5x_{u_{1}}, \quad (q(G)-2)x_{u_{1}}=x_{u}+x_{w}.$$

Thus, $$x_{u}=\frac{q(G)+3}{(q(G)-1)(q(G)-7)}x_{w},\quad\quad x_{u_{1}}=\frac{q(G)-5}{(q(G)-1)(q(G)-7)}x_{w}.$$ Let $G_{12}=G\cup \{v\}$  and $G_{12}^{\prime}=G_{12}-uu_{1}-uu_{2}+u_{1}u_{2}+wv$. It is obvious that $G_{12}^{\prime}=K_{1}\vee ((k-2)P_{2}\cup S_{4}\cup P_{1})$ and $X^{\prime}=(X^{T},0)^{T}$ is an eigenvector of $Q(G_{12})$ corresponding to $q(G_{12})=q(G)$. Thus, we have
$$
\begin{aligned}
q(G_{12}^{\prime})-q(G_{12})&\geq {X^{\prime}}^{T}(Q(G_{12}^{\prime})-Q(G_{12}))X^{\prime}\\
&=(x_{w}+x_{v})^{2}+(x_{u_{1}}+x_{u_{2}})^{2}-(x_{u}+x_{u_{1}})^{2}-(x_{u}+x_{u_{2}})^{2}\\
&=x_{w}^{2}+2x_{u_{1}}^{2}-2x_{u}^{2}-4x_{u}x_{u_{1}}\\
&=\frac{q^{4}-16q^{3}+74q^{2}-136q+141}{(q(G)-1)^{2}(q(G)-7)^{2}}x_{w}^{2}>0,
\end{aligned}
$$
for $q>2k+2+\frac{1}{k}>10, k\geq4$, which implies $q(G)=q(G_{12})<q(K_{1}\vee ((k-2)P_{2}\cup S_{4}\cup P_{1}))$, a contradiction. If there exists a vertex $w^{\prime}$ satisfying $d(w^{\prime})=1$, then $2m=6k+4=\sum_{u\in V(G)}d_{G}(u)\geq 2k+1+(2k-1)\times2=6k-1$, and hence there exists at most a vertex with degree 7 of $G$. We consider these cases in the following.

If $$2m=6k+4=2k+1+2x+3(2k-1-x)),$$ then $x=2k-6$ and $2k-1-x=5$ and hence the degree sequence of $G-w$ is $\{2,2,2,2,2,\underbrace{1,1,\cdots, 1}_{2k-6},0\}$. Thus, $G\in \{K_{1}\vee ((k-4)P_{2}\cup P_{7}\cup P_{1}), K_{1}\vee((k-5)P_{2}\cup P_{6}\cup S_{3}\cup P_{1}),K_{1}\vee((k-5)P_{2}\cup P_{5}\cup P_{4}\cup P_{1}),K_{1}\vee((k-6)P_{2}\cup P_{5}\cup 2S_{3}\cup P_{1}),K_{1}\vee((k-6)P_{2}\cup 2P_{4}\cup S_{3}\cup P_{1}), K_{1}\vee((k-7)P_{2}\cup P_{4}\cup 3S_{3}\cup P_{1}),K_{1}\vee (k-8)P_{2}\cup 5S_{3}\cup P_{1})\}$.

If $$2m=6k+4=2k+1+2x+4(2k-1-x),$$ then $x=\frac{4k-7}{2}$, a contradiction.

If $$2m=6k+4=2k+1+2x+5(2k-1-x),$$ then $x=\frac{6k-8}{3}$, a contradiction.

If $$2m=6k+4=2k+1+2x+6(2k-1-x),$$ then $x=\frac{8k-11}{4}$, a contradiction.

If $$2m=6k+4=2k+1+2x+7(2k-1-x),$$ then $x=2k-2$, and $2k-1-x=1$ and hence the degree sequence of $G-w$ is $\{6, \underbrace{1,1,\cdots, 1}_{2k-2},0\}$. Thus, $G\in \{K_{1}\vee ((k-4)P_{2}\cup S_{7}\cup P_{1})\}$.

If $$2m=6k+4=2k+1+3x+4(2k-1-x),$$ then $x=4k+7<2k-1$, a contradiction.

Note that $\{K_{1}\vee ((k-4)P_{2}\cup P_{7}\cup P_{1}), K_{1}\vee((k-5)P_{2}\cup P_{6}\cup S_{3}\cup P_{1}),K_{1}\vee((k-5)P_{2}\cup P_{5}\cup P_{4}\cup P_{1}),K_{1}\vee((k-6)P_{2}\cup P_{5}\cup 2S_{3}\cup P_{1}),K_{1}\vee((k-6)P_{2}\cup 2P_{4}\cup S_{3}\cup P_{1}), K_{1}\vee((k-7)P_{2}\cup P_{4}\cup 3S_{3}\cup P_{1}),K_{1}\vee (k-8)P_{2}\cup 5S_{3}\cup P_{1}), K_{1}\vee ((k-4)P_{2}\cup S_{7}\cup P_{1})\}\in K_{1}\vee(F\cup P_{1})$, where $F\in \mathcal{F}_{k+2,k-3}$. By Corollary \ref{le:ch-2.10.}, $H\cong (k-4)P_{2}\cup S_{7}\cup P_{1}$ and hence $G\cong K_{1}\vee ((k-4)P_{2}\cup S_{7}\cup P_{1})$. If $G\cong K_{1}\vee ((k-4)P_{2}\cup S_{7}\cup P_{1})$, then let $V(S_{7})=\{u, u_{1},u_{2},u_{3}, u_{4},u_{5},u_{6}\}$ and $u$ be the center vertex of $S_{7}$. Let $X$ be the perron vector of $Q(G)$, we have $x_{u_{1}}=x_{u_{i}}$ for $i\in \{2,3,4,5,6\}$ and $$(q(G)-7)x_{u}=x_{w}+6x_{u_{1}}, \quad (q(G)-2)x_{u_{1}}=x_{u}+x_{w}.$$

Thus, $$x_{u}=\frac{q(G)+4}{(q(G)-1)(q(G)-8)}x_{w},\quad\quad x_{u_{1}}=\frac{q(G)-6}{(q(G)-1)(q(G)-8)}x_{w}.$$ Let $G_{13}=G\cup \{v\}$  and $G_{13}^{\prime}=G_{13}-uu_{1}-uu_{2}-uu_{3}+u_{1}u_{2}+u_{3}v+vw$. It is obvious that $G_{13}^{\prime}=K_{1}\vee ((k-2)P_{2}\cup S_{4}\cup P_{1})$ and $X^{\prime}=(X^{T},0)^{T}$ is an eigenvector of $Q(G_{13})$ corresponding to $q(G_{13})=q(G)$. Thus, we have
$$
\begin{aligned}
q(G_{13}^{\prime})-q(G_{13})&\geq {X^{\prime}}^{T}(Q(G_{13}^{\prime})-Q(G_{13}))X^{\prime}\\
&=(x_{w}+x_{v})^{2}+(x_{u_{1}}+x_{u_{2}})^{2}+(x_{u_{3}}+x_{v})^{2}-(x_{u}+x_{u_{1}})^{2}\\
&-(x_{u}+x_{u_{2}})^{2}-(x_{u}+x_{u_{3}})^{2}\\
&=x_{w}^{2}+2x_{u_{1}}^{2}-3x_{u}^{2}-6x_{u}x_{u_{1}}\\
&=\frac{q^{4}-18q^{3}+90q^{2}-180q+232}{(q(G)-1)^{2}(q(G)-8)^{2}}x_{w}^{2}>0,
\end{aligned}
$$
for $q>2k+2+\frac{1}{k}>12, k\geq5$, which implies $q(G)=q(G_{13})<q(K_{1}\vee ((k-2)P_{2}\cup S_{4}\cup P_{1}))$, a contradiction.

{\bf Subcase 3.1.2.} $H$ contains a cycle $u_{1}u_{2}\cdots u_{s}u_{1}$.

If $\delta(G)\geq2$, then $$2(k+2)=2|E(H)|\geq \sum_{i=1}^{s}d_{H}(u_{i})+(2k-s)\times 1=2k+s,$$ Thus $3\leq s\leq4$. If $s=3$, then there exist a vertex with degree $2$ in $G-w$. Thus, $H\cong (k-3)P_{2}\cup C_{3}\cup S_{3}$ and $G\cong K_{1}\vee ((k-3)P_{2}\cup C_{3}\cup S_{3})$. If $G\cong K_{1}\vee ((k-3)P_{2}\cup C_{3}\cup S_{3})$, then let $V(C_{6})\cup V(S_{3})=\{u_{1},u_{2},u_{3}\}\cup  \{z,z_{1},z_{2}\}$ and $z$ be the center vertex of $S_{3}$. Let $X$ be the perron vector of $Q(G)$, we have $x_{u_{1}}=x_{u_{i}}$ for $i\in \{2,3\}, x_{z_{1}}=x_{z_{2}}$ and $$(q(G)-3)x_{u_{1}}=x_{w}+2x_{u_{1}}.$$Thus, $$x_{u_{1}}=\frac{1}{q(G)-5}x_{w}.$$Let $G_{14}=G\cup \{v\}$ and $G_{14}^{\prime}=G_{14}-u_{1}u_{2}-u_{1}u_{3}+u_{1}z+wv$. It is obvious that $G_{14}^{\prime}=K_{1}\vee ((k-2)P_{2}\cup S_{4}\cup P_{1})$ and $X^{\prime}=(X^{T},0)^{T}$ is an eigenvector of $Q(G_{14})$ corresponding to $q(G_{14})=q(G)$. Thus, we have
$$
\begin{aligned}
q(G_{14}^{\prime})-q(G_{14})&\geq {X^{\prime}}^{T}(Q(G_{14}^{\prime})-Q(G_{14}))X^{\prime}\\
&=(x_{w}+x_{v})^{2}+(x_{u_{1}}+x_{z})^{2}-2(x_{u_{1}}+x_{u_{2}})^{2}\\
&=x_{w}^{2}-7x_{u_{1}}^{2}+x_{z}^{2}+2x_{u_{1}}x_{z}>0,
\end{aligned}
$$
for $q>2k+2+\frac{1}{k}>8, k\geq3$, which implies $q(G)=q(G_{14})<q(K_{1}\vee ((k-2)P_{2}\cup S_{4}\cup P_{1}))$, a contradiction. If $s=4$, then $H\cong (k-2)P_{2}\cup C_{4}$ and hence $G\cong K_{1}\vee ((k-2)P_{2}\cup C_{4})$. If $G\cong K_{1}\vee ((k-2)P_{2}\cup C_{4})$, then let $V(C_{4})=\{u_{1},u_{2},u_{3},u_{4}\}$. Let $X$ be the perron vector of $Q(G)$, we have $x_{u_{1}}=x_{u_{i}}$ for $i\in \{2,3,4\}$ and $$(q(G)-3)x_{u_{1}}=x_{w}+2x_{u_{1}}.$$Thus, $$x_{u_{1}}=\frac{1}{q(G)-5}x_{w}.$$Let $G_{15}=G\cup \{v\}$ and $G_{15}^{\prime}=G_{15}-u_{1}u_{2}-u_{1}u_{4}+u_{1}u_{3}+wv$. It is obvious that $G_{15}^{\prime}=K_{1}\vee ((k-2)P_{2}\cup S_{4}\cup P_{1})$ and $X^{\prime}=(X^{T},0)^{T}$ is an eigenvector of $Q(G_{15})$ corresponding to $q(G_{15})=q(G)$. Thus, we have
$$
\begin{aligned}
q(G_{15}^{\prime})-q(G_{15})&\geq {X^{\prime}}^{T}(Q(G_{15}^{\prime})-Q(G_{15}))X^{\prime}\\
&=(x_{w}+x_{v})^{2}+(x_{u_{1}}+x_{u_{2}})^{2}-2(x_{u_{1}}+x_{u_{2}})^{2}\\
&=x_{w}^{2}-4x_{u_{1}}^{2}>0,
\end{aligned}
$$ for $q>2k+2+\frac{1}{k}>8, k\geq3$, which implies $q(G)=q(G_{15})<q(K_{1}\vee ((k-2)P_{2}\cup S_{4}\cup P_{1}))$, a contradiction.

If there exists a vertex $w^{\prime}$ satisfying $d(w^{\prime})=1$, then $$2(k+2)=2|E(H)|\geq \sum_{i=1}^{s}d_{H}(u_{i})+(2k-1-s)\times 1=2k+s-1,$$ and hence $s\leq5$. We claim that $s=3$. Otherwise, assume that $x_{u_{1}}\geq x_{u_{i}}$ for $u_{i}\in \{u_{2}, \cdots, u_{s}\}$. Let $H^{\prime}=H-u_{3}u_{4}+u_{1}u_{3}$. By Lemma \ref{le:ch-2.2.}, we have $q(K_{1}\vee H)<q(K_{1}\vee H^{\prime})$, a contradiction. Let $U_{1}$ be the component of $H$ that contains the cycle $u_{1}u_{2}u_{3}u_{1}$. We claim that $U_{1}$ is a unicycle graph with at most 4 vertices. As $|E(U_{1})|\geq |V(U_{1})|$ and $$0\leq 2|E(H-U_{1})|-|V(H-U_{1})|=2(k+2-E(U_{1}))-(2k-|V(U_{1})|))=|V(U_{1})|-2|E(U_{1})|+4,$$ we have that $|E(U_{1})|=|V(U_{1})|\in \{3,4\}$ (the case $|E(U_{1})|=4$ and $|V(U_{1})|=3$ is impossible). Thus we have $U_{1}=C_{3}$ or $U_{1}=S_{4}^{+}$. If $U_{1}=C_{3}$, then $H\cong tP_{2}\cup C_{3}\cup P_{1}$ and hence $t=k-2=k-1$, a contradiction. If $U_{1}=S_{4}^{+}$, then $H\cong tP_{2}\cup S_{4}^{+} \cup P_{1}$ and hence $t=k-2=\frac{2k-5}{2}$, a contradiction.

{\bf Subcase 3.1.3.} There exists $v\in V(G)$ such that $v\notin N_{G}(w)$.

If there exists a vertex $v$ with $d_{G}(v)\geq2$, then $$\sum_{u\in N_{G}(w)}d_{G}(u)\leq 2(3k+2)-2k-d(v)\leq 4k+2.$$ Then we have $$q(G)\leq d_{G}(w)+m(w)\leq 2k+\frac{4k+2}{2k}=2k+2+\frac{1}{k},$$ contradicts the inequation \eqref{eq:ch-1}. Thus there exists a vertex $v\notin N_{G}(w)$ with $d_{G}(v)=1$ and let $N(v)=\{u\}$. If $d_{G}(u)\geq 3$, then let $H^{\prime}=H-uv+wv$, By Lemma \ref{le:ch-2.2.}, we have $q(K_{1}\vee H)<q(K_{1}\vee H^{\prime})$, a contradiction. Thus, $ d_{G}(u)=2$ and $N_{G}(u)=\{w, v\}$.  Let $H^{\star}=G-\{w,u,v\}$. Then $2|E(H^{\star})|=6k+4-2k-2-1=4k+1$, a contradiction. $\blacksquare$

\section*{Data availability}

No data was used for the research described in the article.

\section*{Declaration of competing interest}

The authors declare that they have no conflict of interest.

\end{document}